\def\Bbb#1{{\bf #1}}

\def\fnote#1{\footnote}
\def\blacksquare{\hbox{\vrule width 4pt height 4pt depth 0pt}}

\def\cwleftpar#1#2{\leftskip #1 \rightskip #2 plus 1fill}
\def\cwrightpar#1#2{\leftskip #1 plus 1fill \rightskip #2}
\def\cwcenterpar#1#2{\leftskip #1 plus 1fill \rightskip #2 plus 1fill}
\def\cwfullpar#1#2{\leftskip#1\rightskip#2}

\def\cwoutdent#1#2{\llap{\hbox to #1{#2 \hss}}\ignorespaces}
\def\cwparbegin#1#2#3#4#5{
	\ifcase #1 \cwleftpar{#2}{#3}
	\or \cwrightpar{#2}{#3}
	\or \cwcenterpar{#2}{#3}
	\else \cwfullpar{#2}{#3}\fi
	\ifcase #4 \baselineskip = 1.5\baselineskip
	\or \baselineskip = 2\baselineskip
	\or \baselineskip = 3\baselineskip
	\else \baselineskip = 1\baselineskip\fi
	\ifdim #5 > 0in \else \noindent \fi
	\noindent\ignorespaces}
\documentclass{article}
\begin{document}
\advance \vsize by -2\baselineskip
\def\makeheadline{
{\noindent \folio \par
\vskip \baselineskip }}
\vspace{2ex}

\noindent {\Huge Special bases for derivations of\\ tensor algebras}\\[0.4ex]
\noindent {\Large II. The case along paths}
\vspace{3ex}

\noindent Bozhidar Z. Iliev
\fnote{0}{\noindent $^{\hbox{}}$Permanent address:
Laboratory of Mathematical Modeling in Physics,
Institute for Nuclear Research and \mbox{Nuclear} Energy,
Bulgarian Academy of Sciences,
Boul.\ Tzarigradsko chauss\'ee~72, 1784 Sofia, Bulgaria\\
\indent E-mail address: bozho@inrne.bas.bg\\
\indent URL: http://theo.inrne.bas.bg/$^\sim$bozho/}

\vspace{3ex}

{\bf \noindent Published: Communication JINR, E5-92-508, Dubna, 1992}\\[1ex]
\hphantom{\bf Published: }
http://www.arXiv.org e-Print archive No.~math.DG/0304157\\[2ex]

2003 PACS numbers: 02.40.Ma, 02.40.Vh,  04.20.Cv, 04.90.+h\\[2ex]

\noindent
{\small
The \LaTeXe\ source file of this paper was produced by converting a
ChiWriter 3.16 source file into
ChiWriter 4.0 file and then converting the latter file into a
\LaTeX\ 2.09 source file, which was manually edited for correcting numerous
errors and for improving the appearance of the text.  As a result of this
procedure, some errors in the text may exist.
}\\[2ex]

	\begin{abstract}
The existence of local bases in which the components of derivations of
tensor algebras over a differentiable manifold vanish along paths is proved.
The holonomicity of these bases is investigated. The obtained results are
applied to the case of linear connections. Some relations with the
equivalence principle are shown.
	\end{abstract}\vspace{3ex}

 {\bf I. INTRODUCTION}
\nopagebreak
\medskip

The existence of local coordinates in which the components of symmetric
linear connections [1] vanish along a smooth path without selfintersections
is a known classical result $[2, 3]$. In connection with the intensive use of
nonsymmetric linear connections $[1, 2]$ in different physical theories it is
natural these results to be generalized to the case of nonvanishing torsion.

This paper investigates the mentioned problem from the more general viewpoint of arbitrary derivations of the tensor algebra over a differentiable manifold $[1, 2]$. In it the existence of special bases (or coordinates) is proved in which the defined below components of derivations, in particular linear connections, vanish along some path. The question when these bases are holonomic or anholonomic [2] is investigated. Also the connection of these topics with the equivalence principle in general relativity [4] and other gravity theories (see, e.g., [5]) is briefly discussed.

For the explicit mathematical formulation of our problem we shall remind from [6] some facts concerning derivations of the tensor algebra over a manifold.

By proposition 3.3 from $ch. I$ in [1] any derivation $D$ of the tensor algebra over a differentiable manifold $M$ admits a unique representation in the form $D=L_{X}+S$,in which $L_{X}$is the Lie derivative along the vector field $X$ and $S$ is tensor field of type (1,1) considered here as a derivation [1]. Both $X$ and $S$ are uniquely determined by D.

If $S$ is a map from the set of $C^{1}$vector fields into the tensor fields
of type (1,1) and $S:X\mapstochar \to S_{X}$, then the equation
\[
D^{S}_{X}=L_{X}+S_{X}\qquad (1)
\]
 defines a derivation of the tensor algebra over $M$ for any $C^{1}$vector
field $X [1]$.As the map $S$ will hereafter be assumed fixed, such a
derivation, i.e. $D^{S}_{X}$will be called an $S$-derivation along $X$ and
will be denoted for brevity by $D_{X}$. An $S$-derivation is a map $D$ such
that $D:X\mapstochar \to D_{X}$, where $D_{X}$is an $S$-derivation along X.

Let $\{E_{i}, i=1,\ldots  n:=\dim(M)\}$ be a (coordinate or noncoordinate $[2, 3])$ local basis of vector fields in the tangent to $M$ bundle. It is holonomic (resp. anholonomic) if the vectors $E_{1}, \ldots  , E_{n}$commute (resp. don't commute) [2].

The local components $(W_{X})^{i}_{.j}$of $D_{X}$with respect to $\{E_{i}\}$
are defined by the equality
\[
 D_{X}(E_{j})=(W_{X})^{i}_{.j}E_{i}.\qquad (2)
\]
 With respect to $D_{X}$they play the same role as the components of a
linear connection with respect to the connection (see below or $[1, 2])$. It
can easily be computed (see (1) and [1]) that
\[
 (W_{X})^{i}_{.j}:=(S_{X})^{i}_{.j}-E_{j}(X^{i})+C^{i}_{.kj}X^{k}\qquad (3)
\]
in which with $X(f)$ we denote the action of $X=X^{k}E_{k}$on
$C^{1}$scalar functions $f$, i.e. $X(f):=X^{k}E_{k}(f)$, and
$C^{i}_{.kj}$defines the commutators of the basic vectors, i.e.
\[
 [E_{j},E_{k}]=C^{i}_{.jk}E_{i}.\qquad (4)
\]
If we make the change $\{E_{i}\}  \to \{E_{i^\prime }:=A^{i}_{i^\prime
}E_{i}\}$, where $A:=  [A^{i}_{i^\prime }] =: [ A^{i^\prime }_{i}]^{-1}$ is a
nondegenerate matrix function, then from (2) we can see that
$(W_{X})^{i}_{.j}$ transform into
\[
 (W_{X})^{i^\prime }_{..j^\prime }=A^{i^\prime }_{i}A^{j}_{j^\prime
}(W_{X})^{i}_{.j}+A^{i^\prime }_{i}X(A^{i}_{j^\prime }),\qquad (5)
\]
 which, if we introduce the matrices $W_{X}:=  (W_{X})^{i}_{.j}  $ and
$W:= (W_{X})^{i^\prime }_{..j^\prime }  $, will read
\[
W=A^{-1  }W_{X}A+X(A),\qquad (5^\prime )
\]
 where as a first matrix index is understood the superscript and
$X(A):=X^{k}E_{k}(A)=  X^{k}E_{k}(A^{i}_{i^\prime })  $.

If $\nabla $ is a linear connection with local components $\Gamma ^{i}_{.jk}($see, e.g., $[1-3])$, it is fulfilled
\[
 \nabla _{X}(E_{j})=(\Gamma ^{i}_{.jk}X^{k})E_{i}.\qquad (6)
\] Hence, comparing (2) and (6), we see that $D_{X}$is a covariant
differentiation along $X$ iff
\[
 (W_{X})^{i}_{.j}=\Gamma ^{i}_{.jk}X^{k}\qquad (7)
\]
 for some functions $\Gamma ^{i}_{.jk}$.

Let $D$ be an $S$-derivation and $X$ and $Y$ be vector fields.  The torsion
operator $T^{D}$of $D$ is defined as
\[
 T^{D}(X,Y):=D_{X}Y-D_{Y}X-[X,Y],\qquad (8)
\]
 which in a case of linear connection reduces to the corresponding
classical quantity (see $[1-3]$ and below (11)). The $S$-derivation $D$ will
be called torsion free if $T^{D}=0$.

 By means of (1) and (2), we find the local expressions:
\[
(T^{D}(X,Y))^{i}=
(W_{X})^{i}_{.l}Y^{l}-(W_{Y})^{i}_{.l}X^{l}-C^{i}_{.kl}X^{k}Y^{l}.\qquad (9)
\]
 If $D$ is a linear connection $\nabla $, then, due to (7), we have:
\[
(T^{\nabla }(X,Y))^{i}=T^{i}_{.kl}X^{k}Y^{l},\qquad (10)
\]
where $[1-3]$
\[
T^{i}_{.kl}:=-(\Gamma ^{i}_{.kl}-\Gamma ^{i}_{.lk})-C^{i}_{.kl},\qquad (11)
\]
are the components of the torsion tensor of $\nabla $.

The task of this work is the problem to be investigated when along a given
path $\gamma :J  \to M, J$ being a real interval, exist special bases
$\{E_{i^\prime }\}$ in which the components $W$    of an $S$-derivation $D$
along some or all vector fields $X$ vanish. In other words, we are going to
solve $eq. (5^\prime )$ with respect to A under certain conditions, which
will be specified below.

\medskip
\medskip
 {\bf II. Derivations along arbitrary vector fields}

\medskip
In this section we shall investigate the problems for existence and some properties of special bases $\{E_{i^\prime }\}$ in which the components of a given $S$-derivation $D_{X}$along {\it arbitrary vector field} vanish along a path $\gamma :J  \to M$, with $J$ being an ${\Bbb R}$-interval. Let's note that $\{E_{i^\prime }\}$ are supposed to be defined {\it in a neighborhood} of $\gamma (J)$, while the components of $D$ vanish {\it on} $\gamma (J)$.

{\bf Proposition 1.} Along a path $\gamma :J  \to M$ the $S$-derivation is a linear connection if and only if along $\gamma $ exists a basis in which the components of $D$ along every vector field vanish along $\gamma $, i.e. on $\gamma (J)$.

{\bf Remark.} We say that along $\gamma$  $D$ is a linear connection if in
some, and hence in any, basis $\{E_{i}\}$, we have $(cf. (7))$
\[
 W_{X}(\gamma (s))=\Gamma _{k}(\gamma (s))X^{k}(\gamma (s))\qquad (12)
\]
 for some defined on $\gamma (J)$ matrix functions $\Gamma _{k}$and every
vector field $X$, i.e. if (7) is valid for $x\in \gamma (J)$, but for
$x\not\in \gamma (J)$ eq.~(7) may not be true.

{\it Proof.} Let along $\gamma $ the derivation $D$ be a linear connection,
i.e. (12) to be valid. Let at first assume that $\gamma $ is without
selfintersections and that $\gamma (J)$ is contained in only one coordinate
neighborhood $U$ in which some local coordinate basis  $\{E_{i}=\partial /\partial x^{i}\}$ is fixed.

Due to $(5^\prime )$ we have to prove the existence of a matrix
$A=[A^{i}_{i^\prime }]=:[A^{i^\prime }_{i}] ^{-1}$ for which in the basis
$\{E_{i^\prime }=A^{i}_{i^\prime }E_{i}\} W$ the equation $(\gamma (s))=0$ to
be fulfilled for every $X=X^{k}E_{k}$. Substituting (12) into $(5^\prime )$,
we see that the last equation is equivalent to
\[
 \Gamma _{k}(\gamma (s))A(\gamma (s))+E_{k}(A)\mid _{\gamma (s)}=0,
\quad E_{k}=\partial /\partial x^{k}.\qquad (13)
\]
 The general solution of this equation can be constructed
as follows.

Let $V:=J  \times \cdot\times J$, where $J$ is taken $n-1=\dim(M)-1$ times.
Let us fix a one-to-one onto $C^{1}$ map $\eta :J\times V  \to M$ such that
$\eta (\cdot ,\mathbf{ t}_{0})=\gamma $ for some fixed $\mathbf{ t}_{0}\in V$, i.e.
$\eta (s,\mathbf{ t}_{0})=\gamma (s), s\in J$, which is possible iff $\gamma $ is
without selfintersections. In $U \cap \eta (J,V)$ we introduce coordinates
$\{x^{i}\}$ by putting $(x^{1}(\eta (s,\mathbf{ t})),\ldots  ,x^{n}(\eta (s,\mathbf{
t})))=(s,\mathbf{ t}), s\in J, \mathbf{ t}\in V$, which again is possible iff $\gamma
$ is without selfintersections.

If we expand $A(\eta (s,\mathbf{ t}))$ into power series with respect to
 $(\mathbf{  t}-\mathbf{ t}_{0})$, we find the general solution of (13) in the form
 \[
  A(\eta (s,\mathbf{ t}))
=\{{\Bbb I}- \sum_{k=2}^{n}\Gamma _{k}(\gamma (s))
 [x^{k}(\eta (s,\mathbf{ t}))-x^{k}(\eta (s,\mathbf{ t}_{0}))]
 Y(s,s_{0};-\Gamma_{1}\circ \gamma ) \times
 \]
 \[
 \qquad   \times B(s_{0},\mathbf{ t}_{0};\eta )+B_{kl}(s,\mathbf{ t};\eta )
 [x^{k}(\eta (s,\mathbf{ t}))-x^{k}(\eta (s,\mathbf{ t}_{0}))] \times
 \]
 \[
 \qquad \times [x^{l}(\eta (s,\mathbf{ t}))-x^{l}(\eta (s,\mathbf{ t}_{0}))].\qquad
 (14)
 \]
  Here:${\Bbb I}$ is the unit matrix, $s_{0}\in J$ is fixed, $B$ is any
 nondegenerate matrix function of its arguments, $B_{kl}$ and their
 derivatives are bounded matrix functions when $\mathbf{ t}  \to \mathbf{ t}_{0}$ and
 $Y=Y(s,s_{0};Z)$, with $Z$ being a continuous matrix function of $s$, is the
 unique solution of the matrix initial-value problem [7]
\[
 \frac{dY}{ds} =ZY,\quad Y|_{s=s_{0}}={\Bbb I}.\qquad (15)
 \]
So, a matrix A, and consequently a basis $\{E_{i}\}$, with the needed
 property exist.

If $\gamma (J)$ doesn't lie in only one coordinate neighborhood, then by
means of the above described method we can obtain local bases with the looked
for property in different coordinate neighborhoods which form a neighborhood
of $\gamma (J)$. From these local bases we can construct a global basis along
$\gamma $ with the needed property, but generally this basis will not be
continuous in the regions of intersection of two coordinate neighborhoods.
For example, if for some $\gamma $ there doesn't exist coordinate
neighborhood containing $\gamma (J)$ but there are coordinate neighborhoods
$U^\prime $ and $U^{\prime\prime}$ such that
$\gamma (J)\subset (U^\prime ) \cup (U^{\prime\prime})$, then in $U^\prime $
and $U^{\prime\prime}$ there are (see above) bases $\{E_{i^\prime }\}$ and
$\{E_{i^{\prime\prime}}\}$ in which the components of $D_{X}$ vanish for
every $X$ along $\gamma $ and a global basis $\{E^{0}_{i}\}$ in
$(U^\prime )\cup(U^{\prime\prime})$ with the same property may be obtained
by putting $E^{0}_{i}\mid _{x}=\delta ^{i^\prime }_{i}E_{i^\prime }\mid _{x}$
for $x\in U^\prime $ and $E^{0}_{i}\mid _{x}=\delta
^{i^{\prime\prime}}_{i}E_{i^{\prime\prime}}\mid _{x}$for $x\in
U^{\prime\prime}\backslash U^\prime$
(note that $(U^{\prime\prime})\cup(U^\prime )$ is not empty as $\gamma $ is
$C^{1}$ path).

Analogously, if $\gamma $ has selfintersections, then on any "part" of
$\gamma $ without selfintersections exist bases with the needed property.
From these bases can be constructed a global basis with the same property
along $\gamma$. (At the points of selfintersections of $\gamma $ we can
arbitrary fix these bases to be the bases obtained above for some fix part of
$\gamma $ which is without selfintersections.)

Consequently, if along $\gamma$ $D$ is a linear connection, then in a
neighborhood of $\gamma (J)$ a basis $\{E_{i^\prime }\}$ exists in which the
components of $D_{X}$ for every $X$ vanish along $\gamma $.

On the opposite, let us assume the existence of a basis $\{E_{i^\prime }\}$
in which $W$  =0 for every X. Fixing some basis $\{E_{i}\}$ such that
$E_{i^\prime }=A^{i}_{i^\prime }E_{i}$, from $(5^\prime )$ we find
$(W_{X}A+X(A))\mid _{\gamma (s)}=0$ or $W_{X}(\gamma (s))=
=-[(X(A))A^{-1}]\mid _{\gamma (s)}$ which means that (12) is valid for
$\Gamma _{k}(\gamma (s))=-[(E_{k}(A))A^{-1}]\mid _{\gamma (s)}.\blacksquare
$

{\bf Proposition 2.} All bases in which along a path $\gamma $ the
components of an $S$-derivation along every vector field vanish are obtained
along $\gamma $ from one another by linear transformations whose coefficients
are such that the action of the vectors from these bases on them vanish along
$\gamma  ($i.e. on $\gamma (J))$.

{\it Proof.} If $\{E_{i}\}$  and  $\{E_{i^\prime }\}$  are  such  bases,
then  $W$  $(\gamma (s))=W_{X}(\gamma (s))=0$. So, from $(5^\prime )$ follows
$X(A)\mid _{\gamma (s)}=0$ for  every  $X=X^{k}E_{k}$, i.e. $E_{k}(A)\mid
_{\gamma (s)}=0.\blacksquare $

{\bf Proposition 3.} If along a path $\gamma :J  \to M$ for some $S$-derivation $D$ there is a local {\it holonomic} (on $\gamma (J))$ basis in which the components of $D$ along every vector field vanish on $\gamma (J)$, then $D$ is torsion free on $\gamma (J)$. On the opposite, if $D$ is torsion free on $\gamma (J)$ and a continuous basis with the mentioned property along $\gamma $ exists, then all bases with the same property are holonomic along $\gamma $.

{\bf Remark.} In the second part of the proposition we say that the basis must be continuous. This is necessary as any holonomic basis is such.

{\it Proof:} If $\{E_{i^\prime }\}$ is a basis with the mentioned property,
i.e. $W(\gamma (s))=0$ for every $X$ and $s\in J$, then, using (8), we
find
\(
T^{D}(E_{i^\prime },E_{j^\prime })|_{\gamma (s)}
=-[E_{i^\prime },E_{j^\prime }]|_{\gamma (s)}
\)
and consequently $\{E_{i^\prime }\}$ is holonomic at $\gamma (s)$, i.e.
$[E_{i^\prime },E_{j^\prime }]$ $_{\gamma (s)}=0$, iff $0=T^{D}(X,Y)|_{\gamma
(s)}=X^{i^\prime }(\gamma (s))Y^{j^\prime }(\gamma (s))(T^{D}(E_{i^\prime
},E_{j^\prime })|_{\gamma (s)})$ (see proposition 1 and (12)) for every
vector fields $X$ and $Y$, which is equivalent to $T^{D}$ $_{\gamma (J)}=0$.

On the opposite, let $T^{D}|_{\gamma (J)}=0$. We want  to  prove  that
any basis $\{E_{i^\prime }\}$ along $\gamma $ in which $W$  $(\gamma (s))=0$
is holonomic at $\gamma (s), s\in $J. The holonomicity at  $\gamma (s)$
means  $0=[E_{i^\prime },E_{j^\prime }]|_{\gamma (s)}=[-A^{k^\prime
}_{k}(E_{j^\prime }(A^{k}_{i^\prime })-E_{i^\prime }(A^{k}_{j^\prime
}))E_{k^\prime }]|_{\gamma (s)}$. But (see  proposition  1)  the
existence  of $\{E_{i^\prime }\}$ is equivalent to $W_{X}(\gamma (s))=(\Gamma
_{k}X^{k})$  $_{\gamma (s)}$for every X. These two facts, combined with (8),
show that $(\Gamma _{k})^{i}_{.j}=(\Gamma _{j})^{i}_{.k}$. Using this  and
$(\Gamma _{k}A+\partial A/\partial x^{k})|_{\gamma (s)}=0$ (see the proof
of  proposition  1),  we  find
$E_{j^\prime }(A^{k}_{i^\prime })|_{\gamma (s)}
=-A^{j}_{j^\prime }A^{i}_{i^\prime }(\Gamma _{j})^{k}_{.i}$  $_{\gamma
(s)}=E_{i^\prime }(A^{k}_{j^\prime })$  $_{\gamma (s)}$    and    therefore
$[E_{i^\prime },E_{j^\prime }]$  $_{\gamma (s)}=0 ($see above), i.$e
\{E_{i^\prime }\}$  is holonomic on $\gamma (J).\blacksquare $

It can be proved (see below lemma 7) that for any path $\gamma :J  \to M$
every basis $\{E^{\gamma }_{i}\}$ defined only on $\gamma (J)$ can be
extended to a {\it holonomic} basis $\{E^{h}_{i}\}$ defined in a neighborhood
of $\gamma (J)$ and such that $E^{h}_{i}\mid _{\gamma (J)}=E^{\gamma }_{i}$.
In particular, this is true for the restriction $E^{\gamma }_{i^\prime
}=E_{i^\prime }\mid _{\gamma (J)}$ of the considered above special bases
$\{E_{i^\prime }\}$. But in the general case, the extended holonomic bases
$\{E^{h}_{i^\prime }\}$ will not have the special property that
$\{E_{i^\prime }\}$ has.

\medskip
\medskip
 {\bf III. DERIVATIONS ALONG PATHS}

\medskip
	Let $\gamma :J  \to M, J$ being an ${\Bbb R}$-interval, be $a
C^{1}$path and $X$ be a $C^{1} $vector field defined on a neighborhood of
$\gamma (J)$ in such a way that on $\gamma (J)$ it reduces to the tangent to
$\gamma $ vector field, i.e. $X\mid _{\gamma (s)}=$  $(s), s\in $J. The
restriction of an $S$-derivation $D_{X}$on $\gamma (J)$ we shall call
$(S-)$derivation along $\gamma $ and denote it by ${\cal D}^{\gamma }$. Of
course, ${\cal D}^{\gamma }$generally depends on the values of $X$ outside
$\gamma (J)$, but, as this dependence is insignificant for the following, it
will not be written explicitly. So, if $T$ is a $C^{1}$tensor field in a
neighborhood of $\gamma (J)$, then
\[
({\cal D}^{\gamma }T)\mid _{\gamma (s)}:={\cal D}^{\gamma
}_{s}T:=(D_{X}T)\mid _{\gamma (s)}, X\mid _{\gamma (s)}=\dot\gamma(s).\qquad
(16)
\]
 It is easily seen that ${\cal D}^{\gamma }_{s}T$ depends only on the
values of $T\mid _{x}$for $x\in \gamma (J)$, but not on the ones for
$x\not\in \gamma (J)$, and also that ${\cal D}^{\gamma }$is a generalization
of the usual covariant differentiation along curves (see $[2-4]$ or Sect.
$V)$.

When restricted on $\gamma (J)$, the components of $D_{X}$ will be called
components of ${\cal D}^{\gamma }$.

{\bf Proposition 4.} Along any $C^{1}$path $\gamma :J  \to M$ there exists a
basis $\{E_{i^\prime }\}$ in which the components of a given $S$-derivation
${\cal D}^{\gamma }$along $\gamma $ are zeros.

{\it Proof.} Let $\{E_{i}\}$ be some fixed basis in a neighborhood of
 $\gamma (J)$. We have to prove the existence of transformation $E_{i}  \to
 E_{i^\prime }=A^{i}_{i^\prime }E_{i}$such that $W$  $\mid _{\gamma
 (J)}=0$, which, by $(5^\prime )$, is equivalent to the existence of a matrix
 function $A= [A^{i}_{i^\prime }] $ satisfying along $\gamma $ the equation
 $(A^{-1}(W_{X}A+X(A)))\mid _{\gamma (J)}=0, s\in J$, or
\[
\dot\gamma (A)|_{\gamma (s)}\equiv \frac{d A(\gamma(s))}{ds}
=-W_{X}(\gamma (s))A(\gamma (s))\qquad (17)
 \]
  as $X\mid _{\gamma (s)}=\dot\gamma (s)$.
 Along $\gamma $ the general solution of this equation with respect to $A$ is
\[
  \mathbf{ A}(s;\gamma )=Y(s,s_{0};-W_{X}\circ \gamma )B(\gamma ),\qquad (18)
 \]
 where $Y$ is the unique solution of the initial-value problem $(15),
 s_{0}\in J$ is fixed and $B(\gamma )$ is a nondegenerate matrix function of
 $\gamma $.

So, if we take any matrix function A with the property $A(\gamma (s))=\mathbf{
A}(s;\gamma )$ for some $s_{0}$ and $B$ (e.g. using the notation of the proof
of proposition 1, in any coordinate neighborhood in which $\gamma $ is
without selfintersections, we can put $A(\eta (s,\mathbf{
t}))=Y(s,s_{0};-W_{X}\circ \gamma )     B(s_{0},\mathbf{ t}_{0},\mathbf{ t};\gamma )$
for a fixed nondegenerate matrix function $B)$, we see that A carries out the
needed transformation. Hence, the basis $E_{i^\prime }=A^{i}_{i^\prime
}E_{i}$ is the looked for one.\blacksquare

{\bf Proposition 5.} The bases along $\gamma :J  \to M$ in which the
components of ${\cal D}^{\gamma }$ vanish are obtained from one another by
linear transformations whose coefficients on $\gamma (J)$ are constant or may
depend only on $\gamma $.

{\it Proof.} If $\{E_{i}\}$ and $\{E_{i^\prime }\}$ are such bases, then
$W_{X}(\gamma (s))=W$  $(\gamma (s))=0, X\mid _{\gamma (s)}=\dot\gamma(s)$
and from $(5^\prime )$ follows   $(A)\mid _{\gamma (s)}= =dA(\gamma
(s))/ds=0$, i.e.  $A(\gamma (s))$ is a constant or depends only on the map
$\gamma .\blacksquare $

From propositions 4 and 5 we infer that the requirement for the components
of ${\cal D}^{\gamma }$ to vanish along a path $\gamma $ determines the
corresponding special bases with some arbitrariness only on $\gamma (J)$ and
leaves them absolutely arbitrary outside the set $\gamma (J)$. For this
reason we can talk for bases defined only on $\gamma (J)$ in which the
components of ${\cal D}^{\gamma }$ vanish.

{\bf Proposition 6.} Let for a $C^{1}$path $\gamma :J  \to M$ the defined on
$\gamma (J)$ basis $\{E_{i^\prime }\}$ be such that the components of some
$S$-derivation ${\cal D}^{\gamma }$along $\gamma $ vanish. Let $U$ be a
coordinate neighborhood such that on $U\cap (\gamma (J))\neq\emptyset$ the
path $\gamma $ is without selfintersections. Then there is a neighborhood of
$U\cap (\gamma (J))$ in $U$ in which $\{E_{i^\prime }\}$ can be extended to a
coordinate basis, i.e. in this neighborhood exist local coordinates
$\{x^{i^\prime }\}$ such that $E_{i^\prime }\mid _{\gamma (s)}=\partial
/\partial x^{i^\prime }\mid _{\gamma (s)}$.

{\bf Remark.} Said in other words, this proposition means that locally any
spacial on $\gamma (J)$ basis for ${\cal D}^{\gamma }$can be thought
(extended) as a coordinate, and hence holonomic, one (see proposition 5).

{\it Proof.} The proposition is a trivial corollary from the proof of proposition 4 and the following lemma.

{\bf Lemma 7.} Let the path $\gamma :J  \to M$ be without selfintersections
and such that $\gamma (J)$ is contained in some coordinate neighborhood $U$,
i.e. $\gamma (J)\subset U$, and let $\{E_{i^\prime }\}$ be defined on $\gamma
(J)$ continuous basis, i.e. $E_{i^\prime }\mid _{\gamma (s)}$depends
continuously on s. Then there is a neighborhood of $\gamma (J)$ in $U$ in
which coordinates $\{x^{i^\prime }\}$ exists such that $E_{i^\prime }\mid
_{\gamma (s)}=\partial /\partial x^{i^\prime }\mid _{\gamma (s)}$, i.e.
$\{E_{i^\prime }\}$ can be extended in it to a coordinate basis.\blacksquare

{\bf Proof of lemma 7.}
Let us take a one-to-one onto $C^{1}$ map $\eta
:J\times V \to U, V:=J  \times\cdot\times J (n-1$ times), such that $\eta
(\cdot ,\mathbf{ t}_{0})=\gamma $ for some fixed $\mathbf{ t}_{0}\in V$, i.e. $\eta
(s,\mathbf{ t}_{0})=\gamma (s), s\in J$ (cf. the proof of proposition 1). In the
neighborhood $\eta (J,V)\subset U$ we introduce coordinates $\{x^{i}\}$ by
putting $(x^{1}(\eta (s,\mathbf{ t})),\ldots  ,x^{n}(\eta (s,\mathbf{ t})))=(s,\mathbf{
t}), s\in J, \mathbf{ t}\in $V.
Let the nondegenerate matrix
$A=[ A^{i}_{j^\prime}(s;\gamma )]
=: [ A^{i^\prime }_{j}(s;\gamma )] ^{-1}$ defines the decomposition of
$\{E_{i^\prime }\}$ over $\{\partial /\partial x^{i}\}$, i.e.
\[
E_{j^\prime }|_{\gamma (s)}
=A^{j}_{j^\prime }(s;\gamma )
\frac{\partial}{\partial x^j} \Big|_{\gamma (s)  }.\qquad (19)
\]
 If we define the functions $x^{i^\prime }:\eta (J,V)  \to {\Bbb R}$ by
\[
x^{i^\prime }(\eta (s,\mathbf{ t})):=\delta ^{i^\prime }_{i}x^{i}_{0}
+ \int^{s}_{s_{0}}A^{i^\prime }_{1}(u;\gamma )du+A^{i^\prime }_{j}(s;\gamma )
[x^{j}(\eta(s,\mathbf{ t}))-x^{j}(\gamma (s))]+
\]
\[
f^{i^\prime }_{jk}(s,\mathbf{ t};\gamma )[x^{j}(\eta (s,\mathbf{ t}))-x^{j}(\gamma
(s))][x^{k}(\eta (s,\mathbf{ t}))-x^{k}(\gamma (s))]\qquad (20)
\]
 in which
$s_{0}\in J$ and $x_{0}\in V$ are fixed and the functions $f^{i^\prime
}_{jk}$together with their derivatives are bounded when $\mathbf{ t}  \to \mathbf{
t}_{0}$, then, because of $\eta (\cdot ,\mathbf{ t}_{0})=\gamma $, we find
\[
\frac{\partial x^{i'}}{\partial x^j}\Big|_{\gamma(s)}
=
\frac{\partial x^{i'}}{\partial x^j}\Big|_{\eta(s,\mathbf{t}_0)}
=A^{j}_{j^\prime }(s;\gamma).\qquad (21)
\]
As $\det[A^{i}_{j^\prime }(s;\gamma )] \neq 0,\infty $, from
(21) follows that the transformation $\{x^{i}\}  \to \{x^{i^\prime }\}$ is a
diffeomorphism in some lying in $U$ neighborhood of $\gamma (J)$. So, in this
neighborhood $\{x^{i^\prime }\}$ are local coordinates, the basic vectors of
which on $\gamma (J)$ are (see (21) and (19))
\[
\frac{\partial}{\partial x^{j'}}\Big|_{\gamma(s)}
=
\frac{\partial x^i}{\partial x^{j'}}\Big|_{\gamma(s)}
\frac{\partial}{\partial x^{i}}\Big|_{\gamma(s)}
= A^{i}_{j^\prime }(s;\gamma ) \frac{\partial}{\partial x^i} \Big|_{\gamma(s)}
\]

 Hence $\{x^{i^\prime }\}$ are the coordinates we are looking
for.\blacksquare

Lemma 7 has also a separate meaning: according to it any locally continuous
basis defined on $\gamma (J)$ can locally be extended to a {\it holonomic}
basis in a neighborhood of $\gamma (J)$. Evidently, such an extension can be
done and in an anholonomic way. Consequently, the holonomicity problem for a
basis defined only on $\gamma (J)$ depends on the way this basis is extended
in a neighborhood of $\gamma (J)$.

\medskip
\medskip
 {\bf IV. DERIVATIONS ALONG FIXED VECTOR FIELD}

\medskip
Results, analogous to the ones of Sect. II, are true and with respect to $S$-derivations $D_{X}$along a {\it fixed vector field} $X$, i.e. (see Sect. I) for a fixed derivation. As this case is not very interesting from the view-point of its application, we shall present below without proofs only three results about it.

{\bf Proposition 8.} Along a path $\gamma :J  \to M$ the $S$-derivation
$D_{X}$along a fixed vector field $X$ is a linear connection, i.e. (12) is
valid for that fixed $X$, iff along $\gamma $ exists basis $\{E_{i^\prime
}\}$ in which the components of $D_{X}$vanish on $\gamma (J)$.

{\bf Proposition 9.} The bases along $\gamma $ in which the components of
$D_{X}$for a fixed $X$ vanish on $\gamma (J)$ are connected by linear
transformations whose matrices are such that the action of $X$ on them vanish
on $\gamma (J)$.

In a case of $D_{X}$ for a fixed $X$, the analog of proposition 3 is,
generally, not true. But if for $D_{X}, X$ being fixed, is valid (16) on
$\gamma (J)$, then we can construct a class of $S$-derivations $\{^\prime
D\}$ whose components for every $X$ are given by (12). Evidently, for these
derivations the proposition 3 holds, i.e. we have

{\bf Proposition 10.} If along $\gamma $ for $D_{X}, X$ being fixed, is valid (12) and there is a local holonomic (on $\gamma (J))$ basis in which the components of $D_{X}$vanish on $\gamma (J)$, then the above described derivations $^\prime D$ are torsion free on $\gamma (J)$. Vice versa, if $^\prime D$ are torsion free on $\gamma (J)$ and there exists a continuous bases in which the components of $D_{X}$vanish, then between them exist holonomic ones, but generally not all of them are such.

\medskip
\medskip
 {\bf V. THE CASE OF LINEAR CONNECTIONS}

\medskip
Here we shall apply the results from the previous sections to the special case of a linear connection $\nabla $. As the corresponding proofs are almost evident, they will be omitted.

{\bf Corollary 11.} Along every path $\gamma :J  \to M$ for any linear connection $  \nabla $ there exist bases in which the components of $\nabla $ vanish on $\gamma (J)$. These bases are connected with one another in the way described in proposition 2.

{\bf Corollary 12.} One, and hence any, continuous basis along a path $\gamma :J  \to M$ in which the components of a linear connection $\nabla $ vanish on $\gamma (J)$ is holonomic if and only if $\nabla $ is torsion free on $\gamma (J)$.

{\bf Remark.} When $\gamma $ is without selfintersections and $\gamma (J)$ lies in only one coordinate neighborhood, then holonomic bases (coordinates) in which the components of $\nabla $ vanish on $\gamma (J)$ exist if $\nabla $ is {\it torsion free} and vice versa, which is a well known fact; $cf. [1-3, 7]$.

{\bf Corollary 13.} For a torsion free linear connection $\nabla $ along any path $\gamma :J  \to M$ without selfintersections and lying in only one coordinate neighborhood there exist coordinates, or equivalently holonomic bases, in which the components of $\nabla $ vanish on $\gamma (J)$.

{\bf Remark.} This corollary reproduces a classical theorem that can be found, for instance, in [3] and in $[2], ch$. III, \S8, in the latter of which references to original works are given.

{\bf Corollary 14.} If $\frac{D}{ds}|_{\gamma }:=\nabla_{\dot\gamma }$ is the
associated to  $\nabla $  covariant derivative along the $C^{1}$ path
$\gamma :J  \to M$, then on  $\gamma (J)$  exist  bases, obtained  from  one
another by  linear   transformations   whose coefficients are constant or
depend only  on  $\gamma $,  in  which  the components  of  $\nabla _{\cdot
}$ vanish   (on   $\gamma (J))$.   If   $\gamma $   is   without
selfintersections  and  $\gamma (J)$   lies   in   only   one   coordinate neighborhood, then in some neighborhood of $\gamma (J)$ all of these bases could be extended in a holonomic way. These results  generalize  a series of analogous ones, which began with [4]  and concern  linear connections.

\medskip
\medskip
 {\bf VI. CONCLUSION}

\medskip
The above investigation shows that under sufficiently general conditions for derivations of the tensor algebra over a differentiable manifold $M$ along a path $\gamma :J  \to M$ generally anholonomic bases exist in which the derivation's components vanish on $\gamma (J)$. It is important to be noted that when the derivations are along paths, then the corresponding special bases always can be taken as holonomic (or coordinate) ones.

A feature of the considered here case along paths is its independence of the derivations's curvature, which even wasn't introduced here. The cause for this is the one dimensionality of the paths as submanifolds of M. In connection with this it is interesting the same problems to be investigated as here but on arbitrary submanifolds of $M$, which will be done elsewhere.

At the end we want to say a few words for the relation of the obtained in this paper results with the equivalence principle $[4, 5]$, according to which the gravitational field strength, usually identified with the components of some linear connection, is locally transformable to zero by an appropriate choice of the local (called inertial, normal, Riemannian or Lorentz) reference frame (basis or coordinates). So, from mathematical point of view, the equivalence principle states the existence of local bases in
which the corresponding connection's components vanish. The results of this investigation show the strict validity of this statement along any path (on any curve). Hence, we can make the following two conclusions: Firstly, any gravity theory based on space-time with a linear connection is compatible with the equivalence principle along every path, i.e. in it exist (local) inertial frames along paths; these frames are generally anholonomic, but under some, not very restrictive from physical point of view, conditions on the paths (see lemma 7) there exist such holonomic frames of reference. Secondly, in the gravity theories based on linear connections the equivalence principle must not be considered as a principle as it is identically fulfilled because of their mathematical background.

\medskip
 {\bf ACKNOWLEDGEMENTS}

\medskip
\medskip
The author is grateful to Professor Stancho Dimiev (Institute of Mathematics of Bulgarian  Academy of Sciences) and to Dr. Sava Manoff (Institute for Nuclear Research and Nuclear Energy of Bulgarian Academy of Sciences) for the valuable comments and stimulating discussions. He thanks Professor N. A. Chernikov (Joint Institute for Nuclear Research, Dubna, Russia) for the interest of the problems put in this work.

\medskip
\medskip
 REFERENCES

\medskip
1.  Kobayashi S., K. Nomizu, Foundations of Differential Geometry (Interscience Publishers, New York-London, 1963), vol. {\bf I}.
\par
2.  Schouten J. A., Ricci-Calculus: An Introduction to Tensor Analysis and
its Geometrical Applications (Springer Verlag, Berlin-G\"ottingen-Heidelberg,
1954), 2-nd ed.

3.  Rashevskii P. K., Riemannian Geometry and Tensor Analysis (Nauka, Moscow, 1967) (In Russian).\par
4.  C. W. Misner, Thorne K. S., Wheeler J. A., Gravitation (W. H. Freeman and Company, San Francisco, 1973).\par
5.  von der Heyde P., Lett. Nuovo Cimento, vol.{\bf 14}, No.$7, 250, (1975)$.\par
6.  Iliev B. Z., Special bases for derivations of tensor algebras. I. Cases in a neighborhood and at a point, Communication $E5-92-507$, JINR, Dubna, 1992.\par
7.  Hartman Ph., Ordinary Differential Equations (John Wiley \& Sons, New York-London-Sydney, $1964), ch$.IV, \S1.\par
8.  Lovelock D., H. Rund, Tensors, Differential Forms, and Variational Principals (Wiley-Interscience Publication, John Wiley \& Sons, New York-London-Sydney-Toronto, 1975).\par

\end{document}